\documentclass[amscd]{amsart}       
\usepackage{amscd}
\input xypic
\xyoption{all}

\newtheorem{thm}{Theorem}[section]
\newtheorem{dfn}[thm]{Definition}

\newtheorem{cor}[thm]{Corollary}

\newtheorem{lem}[thm]{Lemma}

\newcommand{\Pf}{{\em Proof}. }
\newtheorem{prop}[thm]{Proposition}

\newcommand{\EPf}{\hfill$QED$}

  {\addvspace{\bigskipamount}\noindent{\bf Example\ }}%
  {}

\newenvironment{tightlist}%
{\begin{list}{$\bullet$}{\parsep=0pt \itemsep=0pt
                         \topsep=6pt \partopsep=\parskip}}%
{\end{list}}

\newcounter{ictr}
\newenvironment{ilist}{\begin{list}
                         {(\roman{ictr})}
                         {\usecounter{ictr}
                          \parsep=0pt \itemsep=0pt
                          \topsep=6pt \partopsep=\parskip
                          \setlength{\leftmargin}{0.5truein}}}
                      {\end{list}}

\newcounter{nctr}

\newcommand{\B}{\ensuremath{\mathcal B}}
\newcommand{\C}{\ensuremath{\mathbb C}}

\newcommand{\E}{\ensuremath{\mathcal E}}
\newcommand{\F}{\ensuremath{\mathcal F}}

\renewcommand{\H}{\ensuremath{\mathcal H}}
\newcommand{\K}{\ensuremath{\mathcal K}}

\newcommand{\N}{\ensuremath{\mathbb N}}

\newcommand{\R}{\ensuremath{\mathbb R}}

\newcommand{\Z}{\ensuremath{\mathbb Z}}

\newcommand{\bakideal}{\triangleright}

\newcommand{\cs}{{C^{\star}}}

\newcommand{\DD}{D\!\!\!\!/}

\newcommand{\dvol}{{\operatorname{dvol}}}
\newcommand{\eps}{{\varepsilon}}
\renewcommand{\epsilon}{{\varepsilon}}
\newcommand{\ev}{{\operatorname{ev}}}

\newcommand{\Hom}{{\operatorname{Hom}}}
\newcommand{\id}{{\operatorname{id}}}
\newcommand{\ideal}{\triangleleft}

\newcommand{\Ker}{{\operatorname{Ker}}}
\newcommand{\Mor}{{\operatorname{Mor}}}

\newcommand{\Range}{{\operatorname{Range}}}

\newcommand{\vlim}[1]{\underset{#1}{\varinjlim}}
\newcommand{\xotimes}{\otimes_{X}}

\newcommand{\lacl}{[\![}
\newcommand{\racl}{]\!]_a}
\newcommand{\rxacl}{]\!]_X^a}
\newcommand{\lhcl}{[\![}
\newcommand{\rhcl}{]\!]}

\newcommand{\rxhcl}{]\!]_X}

\newcommand{\RE}{\ensuremath{{\mathcal R}E}}
\newcommand{\RKK}{\ensuremath{{\mathcal R}KK}}

\begin{document}

\title     {Representable $E$-theory for $C_0(X)$-algebras}
\author    {Efton Park}
\address   {Department of Mathematics, Texas Christian University,
            Box 298900, Fort Worth, TX 76129}
\email     {e.park@tcu.edu}
\author    {Jody Trout}
\address   {Department of Mathematics, Dartmouth College, 6188 Bradley Hall,
Hanover, NH 03755}
\email     {jody.trout@dartmouth.edu}

\thanks{The second author was partially supported by NSF Grant DMS-9706767}


\begin{abstract}
Let $X$ be a locally compact space, and let $A$ and $B$ be
$C_0(X)$-algebras. We define the notion of an asymptotic 
$C_0(X)$-morphism from $A$ to $B$ and construct 
representable \linebreak $E$-theory groups $\RE(X;A,B)$. 
These are the universal groups on the
category of separable $C_0(X)$-algebras that are $C_0(X)$-stable,
$C_0(X)$-homotopy-invariant, and half-exact. If $A$ is $\RKK(X)$-nuclear, these groups
are naturally
isomorphic to Kasparov's representable $KK$-theory groups $\RKK(X;A,B)$.
Applications and examples are also discussed.
\end{abstract}

\maketitle

\section{Introduction}

Let $X$ be a locally compact space. A $\cs$-algebra $A$ is called a
$C_0(X)$-algebra if it is equipped with a nondegenerate action
of $C_0(X)$ as commuting multipliers of $A$. These algebras have many important
applications  since they can be represented as the sections of a
noncommutative bundle of $\cs$-algebras over $X$
\cite{kirchberg-wassermann95,ecterhoff-williams97}.
In particular, any continuous-trace $\cs$-algebra with spectrum $X$ is a
$C_0(X)$-algebra \cite{raeburn-williams98}.
There is a substantial amount of literature on $C_0(X)$-algebras; see the
general treatments in \cite{blanchard96,nilsen96}.

In 1988, Kasparov \cite{kasparov88} defined representable $KK$-theory
groups $\RKK(X;A,B)$, which extend his bivariant $KK$-theory groups
$KK(A,B)$ (where $X = \bullet$) and the representable topological
$K$-theory $RK^0(X)$ of Segal \cite{segal70}. These groups (and their
equivariant versions) have found many applications in geometry, topology and
index theory, most notably in
Kasparov's proof of the Novikov Conjecture for closed discrete subgroups
of finite component Lie groups \cite{kasparov88}. However, these groups are
difficult to compute since they
have six-term exact sequences only for split-exact sequences of
$C_0(X)$-algebras \cite{bauval98}.

In 1989, Connes and Higson \cite{connes-higson89} defined $E$-theory
groups $E(A, B)$ based on asymptotic morphisms from $A$ to $B$. An
asymptotic morphism  $\{\phi_t\} : A \to B$  is a one-parameter family of maps
which satisfy the properties of a $*$-homomorphism in the
limit as $t \to \infty$. There is a natural notion of
homotopy for asymptotic morphisms for which a composition product is
well-defined.
If $A$ is separable, the functor $B \mapsto E(A,B)$ is the
universal functor on the category of separable $\cs$-algebras $B$ which
is stable, homotopy-invariant and half-exact. Thus, there are six-term
exact sequences in both variables for any exact sequence of
separable $\cs$-algebras. There is a natural transformation
$KK(A,B) \to E(A,B)$, which is an isomorphism if $A$ is $K$-nuclear.

In this paper, we extend these constructions to the category of
separable $C_0(X)$-algebras by defining a notion of asymptotic
$C_0(X)$-morphism between two $C_0(X)$-algebras $A$ and $B$. There is a
corresponding notion of $C_0(X)$-homotopy for an asymptotic $C_0(X)$-morphism
for which composition is well-defined. This leads us to define
representable $E$-theory groups $\RE(X;A,B)$. If $X$ is second countable then these
groups have a composition product
$$\RE(X;A,B) \otimes \RE(X;B,C) \to \RE(X;A,C),$$
and we prove that these are the universal such groups on the category
of separable $C_0(X)$-algebras that are $C_0(X)$-stable,
$C_0(X)$-homotopy-invariant, and
half-exact. We also prove that there is a natural transformation
$$\RKK(X;A,B) \to \RE(X;A,B),$$
which is an isomorphism if $A$ is $\RKK(X)$-nuclear (in
the sense of Bauval \cite{bauval98}.) If $X = \bullet$ is a point, we recover
the
original $E$-theory groups of Connes-Higson  $E(A, B) \cong \RE(\bullet; A, B)
$.

In addition to these properties, which generalize the basic properties of
$E$-theory, we show that $\RE$-theory also has the following
properties in common with $\RKK$. For $C_0(X)$-algebras $A$, $B$, $C$, $D$,
there is a tensor product operation
$$\RE(X;A,B) \otimes \RE(X;C,D) \to \RE(X;A \xotimes C, B \xotimes D),$$
where $A \xotimes B$ denotes the balanced tensor product over
$C_0(X)$ \cite{ecterhoff-williams97,blanchard95}. If $p : Y \to X$ is a continuous map, the pullback
construction $A \mapsto p^*A$ of Raeburn-Williams \cite{raeburn-williams85}
(which converts a $C_0(X)$-algebra $A$ to a $C_0(Y)$-algebra $p^*A$)
induces a natural transformation of functors
$$p^* : \RE(X;A,B) \to \RE(Y; p^*A, p^*B).$$
In section 2 we review some definitions and constructions involving
$C_0(X)$-algebras that we will need later.  We define $\RE$-theory in
section 3. The category-theoretic tools for comparing $\RE$-theory and
Kasparov's $\RKK$-theory
are developed in section 4. Finally, in section 5 we give some examples and
applications.  For unital $C_0(X)$-algebras, we define the notion of a
``fundamental class'' in $\RE$-theory.  Also, $\RE$-elements are associated to families
$\{D_x\}_{x \in X}$ of elliptic differential operators
parametrized by the space $X$.  We also discuss invariants of
central bimodules, which have found recent applications in noncommutative geometry and physics.


\section{Review of $C_0(X)$-algebras}

Let $X$ be a locally compact Hausdorff topological space.

\begin{dfn}
A $\cs$-algebra $A$ is called a {\bf central Banach $C_0(X)$-module} if $A$ is
equipped with a
$*$-homomorphism $\Phi : C_0(X) \to ZM(A)$ from $C_0(X)$ into the center
of the multiplier algebra of $A$. If $\Phi(C_0(X)) \cdot A$ is dense in $A$,
we say that $\Phi$ is {\bf nondegenerate}, and we call $A$ a
{\bf $C_0(X)$-algebra}. We shall usually write $f \cdot a$
for $\Phi(f)a$.
\end{dfn}

Note that if $A$ is a central Banach $C_0(X)$-module then $A^X =_{\rm def}
\overline{\Phi(C_0(X)) \cdot A}$ is the maximal
subalgebra of $A$ that is a $C_0(X)$-algebra. Compare the following lemma to
Definition 1.5 \cite{kasparov88}, which assumes
$\sigma$-compactness of $X$.

\begin{lem}
The following are equivalent:
\begin{ilist}
\item $A$ is a $C_0(X)$-algebra;
\item There is a $*$-homomorphism $\Phi : C_0(X) \to ZM(A)$ such that for any
approximate unit $\{f_\lambda\}$ in $C_0(X)$
we have $\lim_{\lambda \to \infty} \| \Phi(f_\lambda)a - a \| = 0$ for all $a
\in A$.
\end{ilist}
\end{lem}

Examples of $C_0(X)$-algebras include $\cs$-algebras with Hausdorff spectrum
$X$;
in particular, any continuous-trace $\cs$-algebra is a $C_0(X)$-algebra
\cite{raeburn-williams98}. Also, the
proper algebras occurring in the Baum-Connes Conjecture are examples
\cite{baum-connes-higson94,guentner-higson-trout99}.
In general, a $C_0(X)$-algebra can be realized as the algebra of sections of an
upper-semicontinuous $\cs$-bundle over
$X$ which vanish at infinity
\cite{kirchberg-wassermann95,ecterhoff-williams97}.
See the general discussions in \cite{blanchard96,nilsen96}.

\begin{dfn} Let $A$ and $B$ be $C_0(X)$-algebras (or central Banach
$C_0(X)$-modules). A $*$-homomorphism $\psi : A \to B$ is
called {\bf $C_0(X)$-linear} if $\psi(f \cdot a) = f \cdot \psi(a)$ for all $f
\in C_0(X)$
and $a \in A$. We call such a homomorphism
a {\bf $C_0(X)$-morphism}. We define
{\bf SC*(X)}
to be the category of all separable $C_0(X)$-algebras and
$C_0(X)$-morphisms,
and let $\Hom_X(A,B)$ denote
the set of $C_0(X)$-morphisms from $A$ to $B$.
\end{dfn}

\begin{lem}\label{lem:morphism}
Let $A$ be a $C_0(X)$-algebra and $B$ be a Banach $C_0(X)$-module. If $\psi : A
\to B$
is a $C_0(X)$-morphism, then $\psi(A) \subset B^X$, i.e., $\Hom_X(A,B) = \Hom_X(A, B^X)$.
\end{lem}

\Pf
Since $\psi(f \cdot a) = f \cdot \psi(a) \in B^X$ for all $f \in C_0(X)$ and $a
\in A$, it follows that $\psi(C_0(X) \cdot A) \subset B^X$.
Since the range of $\psi$ is closed and $C_0(X) \cdot A$ is dense in $A$, the
result follows.
\EPf

An easy to prove, but important, result is the following:

\begin{lem}\label{lem:tensor}
If $A$ is a $C_0(X)$-algebra, then under the maximal tensor product,
$A \otimes B$ is a $C_0(X)$-algebra for any
$\cs$-algebra $B$, with the $C_0(X)$-action given by the formula
$f \cdot(a \otimes b) = (f \cdot a) \otimes b$.
\end{lem}

In this paper, we will only use the maximal tensor product,
since it is the appropriate tensor product when working with
asymptotic morphisms.

Let $B$ be any $\cs$-algebra and let $\phi : B \to A$
be a $*$-homomorphism into a
$C_0(X)$-algebra $A$. By the universal
property of the maximal tensor product, the map $f \otimes a \mapsto f \cdot
\phi(a)$ defines a unique
$C_0(X)$-morphism $\phi_X : C_0(X) \otimes B \to A$.

Let $\K$ be the $\cs$-algebra of compact operators on separable Hilbert
space. A
$C_0(X)$-algebra
$A$ is called {\it $C_0(X)$-stable} if $A \cong A \otimes \K$ as
$C_0(X)$-algebras (see the example
following Corollary 6.10 of \cite{raeburn-williams98}
for two $C_0(X)$-algebras with the feature that $A \otimes \K$ and
$B \otimes \K$ are isomorphic as
$\cs$-algebras, but not as $C_0(X)$-algebras.)

\begin{prop}
The category {\bf SC*(X)} is closed under the operations of taking
ideals, quotients, direct sums, suspensions, and $C_0(X)$-stabilizations.
\end{prop}

Let $A$ and $B$ be $C_0(X)$-algebras.
If $i_A : M(A) \to M(A \otimes B)$ and $i_B : M(B) \to M(A \otimes B)$ denote
the natural
inclusions, then $A \otimes B$ has two natural $C_0(X)$-algebra structures.
Typically these $C_0(X)$-algebra structures do not agree.  However,
they can be made to agree by taking a quotient by an appropriate ideal:

\begin{dfn}\label{def:xtensor} \cite{ecterhoff-williams97,blanchard95}
Let $I_X$ be the closed $C_0(X)$-ideal of $A \otimes B$ generated by elementary
tensors
$$\{f \cdot a \otimes b - a \otimes f \cdot b  : f \in C_0(X), a \in A, b \in
B\}.$$
Define $A \xotimes B = (A \otimes B)/I_X$ to be the quotient $\cs$-algebra. It
has a
natural $C_0(X)$-structure given on the images $a \xotimes b$ of elementary
tensors
$a \otimes b$ by
$$f \cdot (a \xotimes b) = (f \cdot a) \xotimes b = a \xotimes (f \cdot b).$$
$A \xotimes B$ is called the (maximal) $C_0(X)$-balanced tensor product of $A$
and $B$.
\end{dfn}

Note that if $A$ is a $C_0(X)$-algebra, then $C_0(X) \xotimes A \cong A$ via
the map $f \xotimes a \mapsto f \cdot a$.

We note that there are other notions of a balanced tensor product,
specifically, the minimal balanced tensor product $\xotimes^m$.
However, this tensor product is not associative (see 3.3.3
\cite{blanchard95}). In addition, the maximal balanced tensor product
enjoys a very useful property that we describe below.
We first observe that if $C$ is a $C_0(X)$-algebra, then the multiplier
algebra $M(C)$ is
canonically a central Banach $C_0(X)$-module (but with a possibly
degenerate $C_0(X)$-action).

\begin{prop} (Compare Prop 2.8 \cite{ecterhoff-williams97})
Let $A$, $B$, and $C$ be $C_0(X)$-algebras. If $\psi : A \to M(C)$ and
$\phi : B \to M(C)$
are commuting  $C_0(X)$-morphisms, then there exists a unique
$C_0(X)$-morphism $$\psi \xotimes \phi : A \xotimes B \to M(C)$$ such that
$(\psi \xotimes \phi)(a \xotimes b) = \psi(a)\phi(b)$.
\end{prop}

Let $p : Y \to X$ be a continuous map of locally compact spaces. The pullback
$\phi(g) = g \circ p$
defines a $*$-homomorphism $\phi : C_0(X) \to C_b(Y) = M(C_0(Y))$ and so gives
$C_0(Y)$ a central
Banach $C_0(X)$-module structure. For a detailed discussion of the following,
see \cite{raeburn-williams85}.

\begin{dfn}
Let $p : Y \to X$ be a continuous map, and let $A$ be a $C_0(X)$-algebra. The
balanced tensor product
$C_0(Y) \xotimes A$ is called the {\rm pullback $\cs$-algebra} and is denoted
$p^*A$. Although it is a $C_0(X)$-algebra, it also has a natural
$C_0(Y)$-action given by the canonical
embedding $C_0(Y) \hookrightarrow M(C_0(Y) \xotimes A)$.
\end{dfn}

If $A = \Gamma_0(E)$ is the $\cs$-algebra of sections of a $\cs$-bundle $E \to
X$ which vanish at infinity,
then $p^*A$ is the $\cs$-algebra of sections of the pullback bundle $p^*E \to
Y$ \cite{raeburn-williams85}.
Note that if $p : Y \to \bullet$ is the map to a point, then $p^*A = C_0(Y)
\otimes A$ for any $\cs$-algebra $A$.

\begin{cor}
Let $p : Y \to X$ be a continuous map. If $\psi : A \to B$ is a
$C_0(X)$-morphism, there is
a canonical $C_0(Y)$-morphism $p^*(\psi) : p^*A \to p^*B$. In fact, the
pullback construction
defines a functor $p^* : $ {\bf SC*(X)} $\to$ {\bf SC*(Y)} that is
both additive and multiplicative.
\end{cor}

Let $B[0,1] = C([0,1], B) \cong B \otimes C[0,1]$ denote the $\cs$-algebra of
continuous functions from the unit interval $[0,1]$ into $B$. It has a
$C_0(X)$-algebra structure by Lemma \ref{lem:tensor}, and the evaluation
maps $\ev_t : B[0,1] \to B$ are $C_0(X)$-morphisms.
Also, the inclusion $i_B : B \to B[0,1]$ by constant functions is a
$C_0(X)$-morphism.  Finally, the flip map
$B[0,1] \to B[0,1] : f(t) \mapsto f(1-t)$ is also a $C_0(X)$-morphism.

\begin{dfn}
Let $A$ and $B$ be $C_0(X)$-algebras. Two $C_0(X)$-morphisms $\psi_0 , \psi_1 :
A \to B$ are
{\bf $C_0(X)$-homotopic} if there is a $C_0(X)$-morphism $\Psi : A \to B[0,1]$
such that
$\ev_i \circ \Psi = \psi_i$ for $i = 0, 1$. This defines an equivalence
relation
on the set of
$C_0(X)$-morphisms from $A$ to $B$. Let $[A,B]_X$ denote the set of
$C_0(X)$-homotopy equivalence
classes of $C_0(X)$-morphisms from $A$ to $B$. We let $[\psi]_X$ denote the
equivalence class
of $\psi : A \to B$, and we define {\bf SH(X)} to be the category of separable
$C_0(X)$-algebras and $C_0(X)$-homotopy classes of $C_0(X)$-morphisms.
\end{dfn}

\section{Representable $E$-theory}

Let $A$ and $B$ be $C_0(X)$-algebras.

\begin{dfn}
An {\rm asymptotic $C_0(X)$-morphism} from $A$ to $B$ is an asymptotic
morphism $\{\psi_t\}_{t \in [1,\infty)} : A \to B$ which is
asymptotically $C_0(X)$-linear, i.e.,
$$\lim_{t \to \infty} \| \psi_t(f \cdot a) - f \cdot \psi_t(a) \| = 0$$
for all $a \in A$ and $f \in C_0(X)$.
\end{dfn}

For a review of the basic properties of asymptotic morphisms of Connes
and Higson, see the papers \cite{connes-higson89,
connes-higson90,guentner-higson-trout99} and the books \cite{connes94,
blackadar98}. A
$C_0(X)$-morphism $\psi : A \to B$ determines a (constant) asymptotic
$C_0(X)$-morphism by the formula $\psi_t = \psi$.

Recall that two asymptotic morphisms $\{\psi_t\}, \{\phi_t\} : A \to B$
are called {\it equivalent} if
$$\lim_{t \to \infty} \| \psi_t(a) - \phi_t(a) \| = 0$$
for all $a \in A$. It is then easy to see that any asymptotic morphism
which is equivalent to an asymptotic $C_0(X)$-morphism is also
asymptotically $C_0(X)$-linear. Let $\lacl \psi_t \rxacl$ denote the
equivalence class of the asymptotic $C_0(X)$-morphism $\{\psi_t\}$, and let
$\lacl A, B \rxacl$ denote the collection of all equivalence classes of
asymptotic $C_0(X)$-morphisms
from $A$ to $B$.

Let $C_b([1,\infty), B)$ denote the continuous bounded functions from the ray
$[1,\infty)$ to $B$. It has an induced structure as a central Banach
$C_0(X)$-module under pointwise multiplication (note, however, that
this action is typically degenerate).
The ideal $C_0([1,\infty), B) \cong C_0([1,\infty)) \otimes B$ of functions
vanishing at infinity is a
$C_0(X)$-algebra by Lemma \ref{lem:tensor}. Therefore, the quotient
$\cs$-algebra
$$B_\infty = C_b([1,\infty), B) / C_0([1,\infty), B)$$
has a canonical central Banach $C_0(X)$-module structure.

\begin{prop}(Compare \cite{connes-higson89})
There is a one-to-one correspondence between equivalence classes of asymptotic
$C_0(X)$-morphisms
$\{\psi_t\} : A \to B$ and $C_0(X)$-morphisms $\Psi : A \to
B_\infty^X$.  In other words,
$$\lacl A, B \rxacl \cong \Hom \sb X(A, B_\infty^X).$$
\end{prop}

\Pf If $\{\psi_t\} : A \to B$ is an asymptotic $C_0(X)$-morphism, we obtain a
$C_0(X)$-morphism $\Psi : A \to B_\infty$ by defining
$$\Psi(a) = q ( \hat{\psi}(a))$$ where $\hat{\psi}(a)(t) = \psi_t(a)$ and
$q : C_b([1,\infty), B) \to B_\infty$ is the quotient map. It depends only 
on the asymptotic equivalence class $\lacl \psi_t \rxacl$. Now invoke Lemma 2.4.
Conversely, given $\Psi \in \Hom_X(A,B^X_\infty)$ we obtain an asymptotic
$C_0(X)$-morphism $\{\psi_t\} : A \to B$ by the composition
$$A \stackrel{\Psi}{\to} B^X_\infty \hookrightarrow B_\infty \stackrel{s}{\to} C_b([1,\infty), B)
\stackrel{\ev_t}{\to} B,$$
where $s$ is any set-theoretic section of $q$. Different choices of $s$ give equivalent asymptotic
$C_0(X)$-morphisms.  Thus, $\lacl A, B \rxacl \cong \Hom_X(A, B_\infty^X)$. \EPf

It follows that the asymptotic equivalence class of an asymptotic
$C_0(X)$-morphism is parametrized
by the choices for the section $s$ which inverts the projection $q :
C_b([1,\infty), B) \to B_\infty$.

\begin{lem}\label{lem:assume}
Let $\{\psi_t\} : A \to B$ be an asymptotic $C_0(X)$-morphism. Then, up to
equivalence,
we may assume that $\{\psi_t\}$ has either one of the following properties (but
not both):
\begin{tightlist}
\item $\{\psi_t\}$ is an equicontinuous family of functions; or
\item Each map $\psi_t$ is $*$-linear.
\end{tightlist}
\end{lem}

\Pf
The first follows from the selection theorem of Bartle and Graves
\cite{bartle-graves52}. The second is realized by using a Hamel basis.
\EPf

\begin{dfn}
Two asymptotic $C_0(X)$-morphisms $\{\psi_t^0\} : A \to B$ and
$\{\psi^1_t\} : A
\to B$ are
{\bf $C_0(X)$-homotopic} if there is an asymptotic $C_0(X)$-morphism
$\{\Psi_t\}
: A \to B[0,1]$
such that $\ev_i \circ \Psi_t = \psi^i_t$
for $i = 0, 1$. We denote this equivalence
by $\{\psi^0_t\} \sim_{hX} \{\psi^1_t\}$, and we let
$\lhcl A, B \rxhcl$ denote the
set of $C_0(X)$-homotopy classes  $\lhcl \psi_t \rxhcl$ of asymptotic
$C_0(X)$-morphisms $\{\psi_t\} : A \to B$.
\end{dfn}

From the discussion in Section 2, it is easy to see that
$C_0(X)$-homotopy defines an equivalence relation on the
set of asymptotic $C_0(X)$-morphisms from $A$ to $B$.  Furthermore,
there is a canonical map
$$[A,B]_X \to \lhcl A, B \rxhcl.$$
In addition, if asymptotic $C_0(X)$-morphisms are equivalent, then they are
actually $C_0(X)$-homotopic via the straight-line homotopy, i.e.,
there is a canonical map
$$\lacl A, B \rxacl \to \lhcl A, B \rxhcl.$$

\begin{lem}\label{lem:reparm}
Let $\{\psi_t\} : A \to B$ be an asymptotic $C_0(X)$-morphism. If $r :
[1,\infty) \to [1, \infty)$ is a continuous
function such that $\lim_{t \to \infty} r(t) = \infty$,
then $\{\psi_{r(t)}\}$ is
an asymptotic $C_0(X)$-morphism
that is $C_0(X)$-homotopic to $\{\psi_t\}$.
\end{lem}

\Pf It is easy to show that $\{\psi_{r(t)}\}$ is an asymptotic morphism.
Let $a \in A$, $f \in C_0(X)$, and $\eps > 0$ be given. By
definition, there is
a $t_0 > 0$ such that
$$\| \psi_t(f \cdot a) - f \cdot \psi_t(a) \| < \eps$$
for all $t \geq t_0$. Now choose $t_1 \geq t_0 > 0$ such that $r(t) \geq t_0$.
Then for
all $t \geq t_1$ we have that
$$\| \psi_{r(t)}(f \cdot a) - f \cdot \psi_{r(t)}(a) \| < \eps.$$
Thus, $\{\psi_{r(t)}\} : A \to B$ is asymptotically $C_0(X)$-linear.
By defining $\{\Psi_t\} : A \to B[0,1]$ via the formula
$$\Psi_t(a)(s) = \psi_{s r(t) + (1-s)t}(a)$$
for all $a \in A$, $t \in [1,\infty)$  and $s \in [0,1]$, we obtain a
$C_0(X)$-homotopy
connecting $\{\psi_t\}$ and $\{\psi_{r(t)}\}$.
\EPf

Let $\{\psi_t\} : A \to B$ be an asymptotic $C_0(X)$-morphism. If $\phi : B \to
C$ is a $C_0(X)$-morphism,
then the composition $\{\phi \circ \psi_t\} : A \to C$ is clearly an asymptotic
$C_0(X)$-morphism. Similarly,
composition on the left with a $C_0(X)$-morphism and $\{\psi_t\}$ produces an
asymptotic $C_0(X)$-morphism.
As for ordinary asymptotic morphisms, composition is a delicate issue, but we
have the following parallel
result for composing asymptotic $C_0(X)$-morphisms. The proof is analogous to
the one given in \cite{connes-higson89},
but with a few extra observations needed for preserving asymptotic
$C_0(X)$-linearity. See also the
expositions in \cite{guentner99, blackadar98}.

Recall that a {\it generating system} $\{A_k\}$ for a separable $\cs$-algebra
$A$ is a countable family of compact
subsets $A_k$ of $A$ such that for all $k$:
\begin{tightlist}
\item $A_k \subset A_{k+1}$;
\item $\cup A_k$ is dense in $A$;
\item each of $A_k + A_k, A_k A_k, A_k^*$ and $\lambda A_k$ for $|\lambda| \leq
1$ is contained in $A_{k+1}$.
\end{tightlist}

\begin{dfn}
Let $X$ be a second countable, locally compact space.  A {\it
$C_0(X)$-generating system} $(\{A_k\}, \{C_k\})$
for a $C_0(X)$-algebra $A$ is a family of compact subsets $\{A_k\}$ of $A$ and
compact subsets $\{C_k\}$ of $C_0(X)$
such that $\{A_k\}$ is a generating system for $A$, $\{C_k\}$ is a generating
system for $C_0(X)$, and
$C_k \cdot A_k \subset A_{k+1}.$
\end{dfn}

The following is then easy to prove, since $C_0(X) \cdot A$ is dense in $A$.

\begin{lem}
Suppose $\{C_k\}$ and $\{A_k\}$ are generating systems for $C_0(X)$ and
$A$, respectively, and let $A_k' = C_k \cdot A_k$.  Then
$(\{A_k'\},\{C_k\})$ is a $C_0(X)$-generating
system for the $C_0(X)$-algebra $A$.
\end{lem}

\begin{thm} (Technical Theorem on Compositions)
Let $X$ be a second countable, locally compact space.
Asymptotic $C_0(X)$-morphisms may be composed
at the level of $C_0(X)$-homotopy; that is, there is an associative map
$$\lhcl A, B \rxhcl \otimes \lhcl B, C \rxhcl \to \lhcl A, C \rxhcl.$$
The composition of $\lhcl \phi_t \rxhcl$ and $\lhcl \psi_t \rxhcl$ is the
$C_0(X)$-homotopy class of an asymptotic $C_0(X)$-morphism $\{\theta_t\}$
constructed from $\{\phi_t\}$ and $\{\psi_t\}$.
\end{thm}

\Pf
Let $(\{A_k\}, \{C_k\})$ be a $C_0(X)$-generating system for the
$C_0(X)$-algebra $A$.
Let $\{\psi_t\} : A \to B$ and $\{\phi_t\} : B \to C$ be asymptotic
$C_0(X)$-morphisms.
We may assume by Lemma \ref{lem:assume} that these asymptotic morphisms are
given by equicontinuous families
of functions.

Since each $A_k$ is compact and $\{\psi_t\}$ is equicontinuous on $A_k$, there
exist $t_k>0$ such that
for all $t \geq t_k$, $a, a' \in A_k$, $f \in C_k$ and $|\lambda| \leq 1$:
\begin{eqnarray*}
\| \psi_t(a a') - \psi_t(a)\psi_t(a') \|       & \leq  1/k \\
\| \psi_t(a + a') - \psi_t(a) - \psi_t(a') \|  & \leq  1/k \\
\| \psi_t(\lambda a) - \lambda \psi_t(a) \|    & \leq 1/k \\
\| \psi_t(a^*) - \phi_t(a)^* \|                & \leq 1/k \\
\| \psi_t(f \cdot a) - f \cdot \psi_t(a) \|    & \leq 1/k \\
\| \psi_t(a) \|   - \|a\|                      & \leq 1/k.
\end{eqnarray*}

Let $(\{B_k\}, \{C_k\})$ be a $C_0(X)$-generating system such that for all $k$,
$$\{\psi_t(a) : t \leq t_k \text{ and } a \in A_{k+10} \} \subset B_k .$$
There then exist $r_k > 0$ such that for all $r \geq r_k$, $b , b' \in
B_{k+10}$, $f \in C_{k+10}$ and $|\lambda| \leq 1$, estimates
similar to the ones above hold for $\{\phi_r\}$. Let $\sigma$ be a continuous
increasing function such that $\sigma(t_k) \geq r_k$.

For each continuous increasing function $\rho \geq \sigma$,
the composition $\theta_t =
\phi_{\rho(t)} \circ \psi_t : A \to C$
defines an asymptotic $C_0(X)$-morphism. We will call $\{\theta_t\}$ a
{\it composition} of $\{\psi_t\}$ and $\{\phi_t\}$.
Different choices of reparametrizations are $C_0(X)$-homotopic through
compositions via Lemma \ref{lem:reparm}.
For different choices of $C_0(X)$-generating systems, it is obviously possible
to choose the function $\sigma$ so that all
estimates needed hold for both $C_0(X)$-generating systems. Hence, for all
$\rho \geq \sigma$ the family $\{\phi_{\rho(t)} \circ \psi_{t}\}$
is a composition for both $C_0(X)$-generating systems.

Consequently, this construction passes to a well-defined pairing
$$\lacl A, B  \rxacl \times \lacl B, C \rxacl \to \lhcl A, C \rxhcl.$$
By taking compositions of $C_0(X)$-homotopies, it passes to a well-defined
pairing on $C_0(X)$-homotopy classes, as desired.
\EPf

Define {\bf AM(X)} to be the category of separable $C_0(X)$-algebras and
$C_0(X)$-homotopy classes of asymptotic
$C_0(X)$-morphisms, i.e., $\Mor(A,B) = \lhcl A, B \rxhcl$. There are obvious
functors {\bf SC*(X)} $\to$ {\bf SH(X)} $\to$ {\bf AM(X)}.

Let $\{\psi_t\} : A \to B$ be an asymptotic $C_0(X)$-morphism. By Lemma
\ref{lem:assume}, we may assume, up
to equivalence, that the family $\{\psi_t\}$ is equicontinuous. Let $SA =
C_0(\R, A) \cong C_0(\R) \otimes A$
denote the suspension of $A$, which is a $C_0(X)$-algebra by Lemma
\ref{lem:tensor}. It follows that the
formula
$$ g \mapsto \phi_t \circ g$$
then defines an asymptotic morphism $\{S\psi_t\} : SA \to SB$ called a
{\it suspension} of $\{\psi_t\}$. It is easily
seen to be asymptotically $C_0(X)$-linear. By applying this construction to
$C_0(X)$-homotopies we obtain
the following result.

\begin{prop}\label{prop:suspend}
Let $A$ and $B$ be $C_0(X)$-algebras. There are well-defined maps
$$\aligned S : \lacl A, B \rxacl & \to \lacl SA, SB \rxacl \\
           S : \lhcl A, B \rxhcl & \to \lhcl SA, SB \rxhcl \endaligned$$
which are induced by sending the class of $\{\psi_t\}$ to the class of
$\{S\psi_t\}$.
\end{prop}

Consider the balanced tensor product $A \xotimes B$ of $A$ and $B$ over
$C_0(X)$. We have the
following generalization of the tensor product for regular asymptotic
morphisms,
which requires
an extra bit of proof.

\begin{prop}\label{prop:tensorX}
Let $\{\psi_t\} : A \to B$ and $\{\phi_t\} : C \to D$ be asymptotic
$C_0(X)$-morphisms.
The tensor product over $X$ is defined (up to equivalence) by the formula
$$\{\psi_t \xotimes \phi_t\} : A \xotimes C \to B \xotimes D: a \xotimes b
\mapsto \psi_t(a) \xotimes \phi_t(b).$$
It is well-defined on equivalence classes and $C_0(X)$-homotopy classes, i.e.,
there are well-defined maps
$$\aligned
\lacl A, B \rxacl \otimes \lacl C, D \rxacl &\to \lacl A \xotimes C, B \xotimes
D \rxacl \\
\lhcl A, B \rxhcl \otimes \lhcl C, D \rxhcl &\to \lhcl A \xotimes C, B \xotimes
D \rxhcl.
\endaligned$$
\end{prop}

\Pf
By Lemma 3.1 of \cite{connes-higson89}, the map $a \otimes b \mapsto \psi_t(a)
\otimes \phi_t(b)$
determines an asymptotic morphism $\{\psi_t \otimes \phi_t\} : A \otimes C
\to B
\otimes D$
on maximal tensor products which is well-defined on equivalence (and homotopy)
classes.
We may assume that each of these maps is $*$-linear up to equivalence by Lemma
\ref{lem:assume}.

Let $I_X$ and $J_X$ be the balancing ideals of $A \otimes C$ and $B \otimes D$,
respectively (see Definition 2.7). Let
$q : B \otimes D \to B \xotimes D = (B\otimes D)/J_X$
denote the quotient map. A simple $\eps/3$-argument then shows that for any
generator
$x = (f \cdot a) \otimes b - a \otimes (f \cdot b)$ of $I_X$ we have
$$\lim_{t \to \infty} \|q( (\psi_t \otimes \phi_t) (x))\| = 0.$$
By \cite{guentner99}, the family $\{\psi_t \otimes \phi_t\}$ determines an
equivalence class of a
{\it relative} asymptotic morphism
$$\{\psi_t \otimes \phi_t \} : A \otimes C \bakideal I_X \to B \otimes D
\bakideal J_X.$$
By Lemma 3.8 in \cite{guentner99} there is an induced asymptotic morphism
on the
quotients
$$\overline{\{\psi_t \otimes \phi_t \}} : (A \otimes C)/I_X \to (B \otimes
D)/J_X$$
which is well-defined up to equivalence and is given by the formula
$a \xotimes b \mapsto \psi_t(a) \xotimes \phi_t(b)$. We shall call this
asymptotic
morphism the {\it tensor product morphism over} $X$ of
$\{\psi_t\}$ and $\{\phi_t\}$ and denote
it by $\{\psi_t \xotimes \phi_t\} : A \xotimes C \to B \xotimes D$.
\EPf

Let $I$ be an ideal in a separable $C_0(X)$-algebra $A$. Then $I$ has a
canonical $C_0(X)$-algebra
structure inherited from $A$. Let $\{u_t\}$ be a quasicentral approximate unit
for $I \ideal A$ \cite{arveson77}.
Let $s : A/I \to A$ be any set-theoretic section of the $C_0(X)$-linear
projection $p : A \to A/I$.
Connes and Higson \cite{connes-higson89} construct an asymptotic morphism
$$\{\alpha_t\} : A/I(0,1) \to I,$$
which on elementary tensors is induced by the formula
\begin{eqnarray}\label{fmla:ext}
\alpha_t : g \otimes \bar{a} \mapsto g(u_t)s(\bar{a})
\end{eqnarray}
for any $g \in C_0(0,1)$ and $\bar{a} = p(a) \in A/I$.
The asymptotic equivalence class $\lacl \alpha_t \racl$ is independent of the
choice of the section $s$. In addition, the (full) homotopy class
$\lhcl \alpha_t \rhcl$ is independent
of the choice of the quasicentral approximate unit
\cite{connes-higson89,blackadar98}.

\begin{lem}\label{lem:amext}
The asymptotic morphism $\{\alpha_t\} : A/I(0,1) \to I$ is asymptotically
$C_0(X)$-linear, and the $C_0(X)$-homotopy class is independent of
the choice of quasicentral approximate unit.
\end{lem}

\Pf
Take $\bar{a} \in A/I$ and $f \in C_0(X)$. By Lemma \ref{lem:tensor},
$A/I(0,1)
= C_0(0,1) \otimes A/I$
has a $C_0(X)$-algebra structure induced from that of $A/I$. For any $\bar{a}
\in A/I$ we have
$f \cdot s(\bar{a}) - s(f \cdot \bar{a}) \in I$, since the projection $p$ is
a $C_0(X)$-morphism.
Therefore, for any $g \in C_0(0,1)$, $f \in C_0(X)$ and $\bar{a} \in A/I$ we
compute that
$$\lim_{t \to \infty} \| \alpha_t( f \cdot (g \otimes \bar{a})) - f \cdot
\alpha_t(g \otimes \bar{a})\|
= \lim_{t \to \infty} \| g(u_t)(f \cdot s(\bar{a}) - s(f \cdot \bar{a}))\|
= 0,$$
since $\{u_t\}$ is quasicentral for $I \ideal A$. It follows that the
associated
$*$-homomorphism $\alpha : A/I(0,1) \to I_\infty$
is $C_0(X)$-linear.

Let $\{u_t\}$ and $\{u_t'\}$ be two quasicentral approximate units for
$I \ideal A$.  Then
the convex combination $u_t^\lambda = \lambda u_t + (1-\lambda)u_t'$ defines a
continuous family for $\lambda \in [0,1]$ of quasicentral
approximate units connecting  $\{u_t\}$ and $\{u_t'\}$. Replacing $\{u_t\}$ in
formula \ref{fmla:ext}
with $\{u_t^\lambda\}$, we obtain a $C_0(X)$-homotopy by the above argument.
\EPf

Let $\psi : A \to B$ be a $C_0(X)$-morphism. The {\it mapping cone} of
$\psi$ is
the $\cs$-algebra $C_\psi$
of $A \oplus B[0,1)$ consisting of all elements $(a, G)$ such that $\psi(a) =
G(0)$. It has an induced
$C_0(X)$-algebra structure given by $f \cdot (a, G) =
(f \cdot a, f\cdot G)$, and there is a commutative diagram
of $C_0(X)$-algebras
$$\CD
C_\psi @>{q}>> B[0,1) \\
@V{p}VV         @VV{\ev_0}V   \\
A      @>>{\psi}> B.
\endCD$$
The $C_0(X)$-morphisms $p$ and $q$ are the projections onto the respective
factors. There is also
a $C_0(X)$-inclusion $j : SB \cong B(0,1) \to C_\psi$ given by $j : F \mapsto
(0,F)$.

The proofs of the following two results are identical to the regular
$E$-theoretic case; one need only
see that the relevant constructions are all asymptotically
$C_0(X)$-linear.

\begin{lem}
Let $\psi : A \to B$ be a $C_0(X)$-morphism. Then for any $C_0(X)$-algebra $D$,
the following sequence is exact:
$$\lhcl D, C_\psi \rxhcl \stackrel{p_*}{\to} \lhcl D, A \rxhcl
\stackrel{\psi_*}{\to} \lhcl D, B \rxhcl.$$
\end{lem}

\Pf
Same proof as Proposition 5 of \cite{dadarlat94} and Proposition 3.7 of
\cite{rosenberg82}.
\EPf

\begin{lem}
Let $\psi : A \to B$ be a $C_0(X)$-morphism. For any $C_0(X)$-algebra $D$,
consider the
following diagram:
$$\diagram
&  \lhcl B, D \rxhcl \dto^{S} \rto^{\psi^*} &
\lhcl A, D \rxhcl \\
\lhcl C_\psi , SD \rxhcl \rto^{j^*} & \lhcl SB, SD \rxhcl          &
\enddiagram$$
Then $S(\Ker(\psi^*)) \subset \Range(j^*).$
\end{lem}

\Pf
Same proof as Proposition 8 of \cite{dadarlat94} and Proposition 3.10 of
\cite{rosenberg82}.
\EPf

Let $0 \to J \to A \stackrel{\psi}{\to} B \to 0$ be a short exact sequence of
$C_0(X)$-algebras, and consider the
$C_0(X)$-morphism $i : J \to C_\psi$ given by $j : a \to (a, 0)$.
There is also
a $C_0(X)$-extension $$0 \to J(0,1) \to A[0,1) \to C_\psi \to 0,$$
where the first map is the obvious inclusion and the map $A[0,1) \to C_\psi$ is
given by the
formula $F \mapsto (F(0), \psi \circ F)$.

\begin{prop}
The suspension $Si : SJ \to SC_\psi$ is an isomorphism in the category {\bf
AM(X)}.  The inverse is
the $C_0(X)$-homotopy class of the asymptotic $C_0(X)$-morphism
$\{\alpha_t\}$ associated to the $C_0(X)$-extension
$$0 \to J(0,1) \to A[0,1) \to C_\psi \to 0$$
from Lemma \ref{lem:amext}. That is,
$\lhcl Si \circ \alpha_t\rxhcl = \lhcl \id_{SC_\psi} \rxhcl$ and  $\lhcl
\alpha_t \circ Si\rxhcl = \lhcl \id_{SJ} \rxhcl$.
\end{prop}

\Pf
Same as the proof of Theorem 13 in \cite{dadarlat94}.
\EPf

\begin{thm} Let $A$ and $B$ be $C_0(X)$-algebras.
\begin{ilist}
\item If $B \cong B \otimes \K$ is $C_0(X)$-stable, then $\lhcl A, B \rxhcl$ is
an abelian monoid with
addition given by direct sum:
$$\lhcl \psi_t \rxhcl + \lhcl \phi_t \rxhcl = \lhcl \psi_t \oplus \phi_t
\rxhcl,$$
using fixed $C_0(X)$-isomorphisms $B \cong B \otimes M_2 \cong M_2(B)$. The
class
of the zero morphism
serves as the identity.
\item $\lhcl A, SB \rxhcl$ is a group under (not necessarily abelian) ``loop
composition'':
given $\{\psi_t\}, \{\phi_t\} : A \to SB$, define $\{\psi_t \cdot \phi_t\} :
A \to SB$ via
$$
(\psi_t \cdot \phi_t)(a)(s) = \left\{
\begin{array}{ll}
   \psi_t(a)(2s), & \mbox{if $0 \leq s \leq 1/2$} \\
   \phi_t(a)(2s-1), & \mbox{if $1/2 \leq s \leq 1.$}
\end{array} \right.$$
\item $\lhcl A, S(B \otimes \K) \rxhcl$ is an abelian group and the
two operations decribed in (i) and (ii) coincide.
\end{ilist}
\end{thm}

\begin{dfn}
Let $A$ and $B$ be $C_0(X)$-algebras. The {\bf representable $E$-theory group}
$\RE(X;A,B)$ is defined
as the direct limit
$$\RE(X;A,B) = \vlim{n \in \N} \ \lhcl S^n(A \otimes \K), S^n(B \otimes \K)
\rxhcl,$$
where the direct limit is taken over the suspension operation of Proposition
\ref{prop:suspend}.
\end{dfn}

We now assume that $X$ is a second countable locally compact space.

\begin{thm}\label{thm:REtheory}
Let $A$ and $B$ be $C_0(X)$-algebras. The representable $E$-theory group
$\RE(X; A,B)$ is a $C_0(X)$-stable,
$C_0(X)$-homotopy invariant and half-exact bifunctor from {\bf SC*(X)} to the
category of abelian groups that
satisfies the following properties:
\begin{ilist}
\item There exists a group homomorphism
$\RE(X; A, B) \otimes \RE(X; B, C) \to \RE(X;A, C)$;
\item There exists a group homomorphism
$$\RE(X; A, B)  \otimes \RE(X; C, D) \to \RE(X; A \xotimes C, B \xotimes D);$$
\item $\RE(X; A,B)$ is an $\RE(X; C_0(X), C_0(X))$-module;
\item $\RE(X;A, B) \cong \RE(X; SA, SB)$;
\item A continuous map $p : Y \to X$ induces a group homomorphism
$p^* : \RE(X; A, B) \to \RE(Y;p^*A, p^*B)$.
\end{ilist}
\end{thm}

\Pf $\RE(X;A,B)$ is $C_0(X)$-stable and $C_0(X)$-homotopy invariant by definition. Half-exactness
follows from Lemmas 3.11, 3.12, 3.13 and Proposition 3.14. Property (i) follows from Theorem 3.8.
Property (ii) follows from Theorem 3.10. The third property follows from (ii) and the 
natural isomorphisms $A \xotimes C_0(X) \cong A$ and $B \xotimes C_0(X) \cong B$. Property (iv)
follows by definition. The last property follows from Definition 2.9 of the pullback 
$p^*A = C_0(Y) \xotimes A$ and the proof of Proposition 3.10.
\EPf

\section{Categorical aspects of $\RKK$-theory and $\RE$-theory}

In this section, we prove universality properties of $\RKK$-theory
and $\RE$-theory; our proofs are largely based on the work of
\cite{higson87} and \cite{blackadar98}.  We begin by recalling the
definition of $\RKK$-theory:

\begin{dfn}\label{def:RKKcycle} Let $A$ and $B$ be $C_0(X)$ algebras.  A
cycle in
$\RKK(X; A, B)$ is a triple $(\E, \phi, F)$, where
\begin{tightlist}
\item $\E$ is a countably generated $\Bbb Z_2$-graded Hilbert
$B$-module;
\item $\phi : A  \rightarrow \B(\E)$ is a $*$-homomorphism;
\item $F \in \B(\E)$ has degree one, and $[F,\phi(a)]$,
$(F - F^*)\phi(a)$, and $(F^2 - 1)\phi(a)$ are all in
$\K(\E)$ for all $a$ in $A$;
\item $\phi(fa)eb = \phi(a)e(fb)$ for all $a \in A$, $b \in B$, $e \in
\E$, and $f \in C_0(X)$.
\end{tightlist}
The abelian group $\RKK(X; A, B)$ is formed by identifying these cycles
under the equivalence relations of unitary equivalence and homotopy.
\end{dfn}

By the same arguments that are used in the
\lq\lq quasihomomorphism picture\rq\rq\  of $KK$-theory
(\cite{blackadar98}, Section 17.6), we may represent
elements of $\RKK(X; A, B)$ by pairs $(\phi^+, \phi^-)$ of
homomorphisms from $A$ to $M(B \otimes\K)$ such that
\begin{tightlist}
\item $\phi^\pm(fa)(b\otimes k) = \phi^\pm(a)(fb\otimes k)$
for all $a \in A$, $b \in B$, $k \in \K$, and $f \in C_0(X)$;
\item $\phi^+(a) - \phi^-(a)$ is in $B \otimes\K$ for each
$a$ in $A$.
\end{tightlist}
We shall call such pairs {\it $C_0(X)$-quasihomomorphisms}.
A $C_0(X)$-quasihomomorphism
is {\it degenerate} if $\phi^+$ and $\phi^-$ are equal, and the
appropriate equivalence relations on
$C_0(X)$-quasihomomorphisms $(\phi^+, \phi^-)$ are
addition of degenerate pairs, conjugation of $\phi^+$ and
$\phi^-$ by the same unitary, and homotopies given by paths
$(\phi^+_t, \phi^-_t)$, where $(\phi^+_t, \phi^-_t)$
is a $C_0(X)$-quasihomomorphism for each $t$ in $[0,1]$ and $a \mapsto \phi^\pm_t(a)$ is continuous
for each $a \in A$.

Let
$$\xymatrix@1{0\ar[r] & B\otimes\K \ar[r]^j & D
\ar[r]^q  & A\ar@/^/[l]^s\ar[r] & 0}$$
be a split exact sequence of $C_0(X)$-algebras, and
let $\iota$ be the
canonical inclusion of $D$ into $M(B \otimes \K)$.  Then $\iota$ is a
$C_0(X)$-homomorphism, and
$$\left(H_B \oplus H_B^{op}, \iota \oplus (\iota\circ s \circ q),
\left(\begin{array}{ccc} 0 & 1 \\ 1 & 0 \end{array}\right)\right)$$
determines an element $\pi_s$ of $\RKK(X; A, B)$ that is called a
{\it splitting morphism}.  Furthermore, every element
of $\RKK(X; A, B)$ is the pullback of a splitting morphism via a
$C_0(X)$-homomorphism:

\begin{prop}(Compare Proposition 17.8.3, \cite{blackadar98})
Let $(\phi^+, \phi^-)$ be a $C_0(X)$-quasi\-homo\-morphism
from $A$ to  $B$.  Then there exists a split exact sequence
$$\xymatrix@1{0\ar[r] & B\otimes\K \ar[r]^j & D
\ar[r]^q  & A\ar@/^/[l]^s\ar[r] & 0}$$
of $C_0(X)$-algebras
and a $C_0(X)$-homomorphism
$f: A \rightarrow D$ such that $f^*(\pi_s) = [\phi^+, \phi^-]$.
\end{prop}

\Pf  Let $D = \{(a, \phi^+(a)+ \beta) :  a \in A, \beta \in B \otimes \K\}$,
and let $q$ be projection onto the
first factor.  The kernel of $q$ is the ideal
$\{(0, \beta) :  \beta \in B \otimes \K\}$, which is of course isomorphic to
$B \otimes \K$.  The map $s: A \rightarrow D$ defined by
$s(a) = (a, \phi^-(a))$ is obviously a $C_0(X)$-splitting for $q$,
and thus we have $\pi_s$ as an element of $\RKK(X; D, B)$.
Define $f: A \rightarrow D$ as $f(a) = (a,\phi^+(a))$; an easy
computation shows that $f^*(\pi_s)  = [\phi^+, \phi^-]$
in $\RKK(X; A, B)$.
\EPf

\begin{thm}\label{thm:RKKcat}
 Let {\bf RKK} denote the category of separable
$C_0(X)$-algebras, with
{\bf RKK} classes serving as morphisms, and let $C$ be the obvious functor
from {\bf SC*(X)} to {\bf RKK}.
 Let {\bf A} be an additive category and let
 $F:${\bf SC*(X)} $\rightarrow$  {\bf A} be a covariant
 functor satisfying the
 following:
 \begin{tightlist}
 \item $F$ is a homotopy functor;
 \item $F$ is stable; i.e., if $e$ is a rank one projection, then the map
  $b \mapsto e \otimes b$ induces an invertible morphism
  $Fe: F(B) \rightarrow F(B \otimes\K)$;
 \item if $\xymatrix@1{0\ar[r] & J \ar[r]^j & D
\ar[r]^q  & A\ar@/^/[l]^s\ar[r] & 0}$  is a split short exact
sequence of $C_0(X)$-algebras, then $F(D)$ is isomorphic to the
direct sum of $F(J)$ and $F(A)$ via the morphisms $j_*$ and
$s_*$.
\end{tightlist}
Then there is a unique functor $\widehat F:$ {\bf RKK}
$\rightarrow$ {\bf A} such that ${\widehat F}\circ C = F$.
\end{thm}

\Pf Since {\bf RKK} and {\bf SC*(X)} consist of precisely the same
objects, the only question is how $\widehat F$ acts on morphisms.
If $x$ is a morphism arising from a $C_0(X)$-homomorphism, then
$\widehat F(x) = F(x)$.  On the other hand,
suppose that $x$ comes from a splitting
morphism represented by the split short exact sequence
$\xymatrix@1{0\ar[r] & B\otimes\K \ar[r]^j & D
\ar[r]^q  & A\ar@/^/[l]^s\ar[r] & 0}$. Then
${\widehat F}(x) \in$
{\bf A}$\left(F(D),F(B)\right)$ is given by the following composition:
$$\xymatrix@1{F(D) \ar[r] ^-{(Fj,Fs)^{-1}} & F(B\otimes\K)
\oplus F(A) \ar[r]^-{p_1} & F(B\otimes\K) \ar[r]^-{Fe^{-1}} & F(B)}.$$
The fact that $F(x)$ is unique follows from the fact that both
$(Fj,Fs)$ and $Fe$ are invertible, and the proof that $\widehat F$ is
well-defined is the same as that in \cite{blackadar98}, Theorem 22.2.1.
\EPf

Next we consider the universality  property of $\RE$-theory.
The proofs of this are almost exactly the same as
the ones given for $E$-theory in
Section 25.6 of \cite{blackadar98}, so we shall only outline the
argument here.
Let $\{\psi_t\}_{t \in [1,\infty)} : A \to B$ be an asymptotic
$C_0(X)$-morphism, let $\Psi: A \to B_\infty$ be the map defined in
Proposition 3.2, and let $q$ be the quotient map from
$C_b([1,\infty), B) $ onto $B_\infty$.  Set
$D = (q^{-1}(\Psi(A))$.  Then we have a short exact sequence
$0 \rightarrow C_0([1,\infty), B) \rightarrow D
\stackrel{\pi}{\rightarrow} A \rightarrow 0$.  Since
$C_0([1,\infty), B)$ is contractible, $\pi$ induces an isomorphism
$\lacl \pi \rxacl$ in $\RE(D,A)$.

\begin{lem}(Compare Proposition 25.6.2, \cite{blackadar98})
 $\lacl \ev_1 \rxacl \circ
 \left({\lacl \pi \rxacl}\right)^{-1}=
 \lacl \psi \rxacl$ in $\RE(X;A,B)$.
\end{lem}

\Pf  For each $t \geq 1$, the asymptotic
$C_0(X)$-morphism $\{\psi_t\circ\pi\}_{t \in [1,\infty)}$ from $D$
to $B$ is equivalent to $\{\ev_t\circ\pi\}_{t \in [1,\infty)}$, which
in turn is homotopic to the constant morphism
$\{\ev_1\circ\pi\}_{t \in [1,\infty)}$.  Therefore
$\lacl \psi\rxacl \circ \lacl \pi \rxacl = \lacl \ev_1\rxacl$ in
$\RE(X;A,B)$, whence the desired result follows.
\EPf

We also need the following lemma, whose proof is easily adapted from
the proof of the \lq\lq classical\rq\rq\  result:

\begin{lem}(\cite{cuntz84}) Let {\bf A} be an additive category and let
 $F:${\bf SC*(X)} $\rightarrow$  {\bf A} be a covariant
 functor satisfying the
 following:
 \begin{tightlist}
 \item $F$ is a homotopy functor;
 \item $F$ stable;
 \item $F$ is half-exact.
\end{tightlist}
Then $F$ satisfies Bott periodicity.
\end{lem}

\begin{thm}\label{thm:REcat}(Compare Theorem 25.6.1, \cite{blackadar98})
 Let {\bf RE} denote the category of separable
$C_0(X)$-algebras, with
{\bf RE} classes serving as morphisms, and let $C$ be the obvious functor
from {\bf SC*(X)} to {\bf RE}.
 Let {\bf A} be an additive category and let
 $F:${\bf SC*(X)} $\rightarrow$  {\bf A} be a covariant
 functor satisfying the
 following:
 \begin{tightlist}
 \item $F$ is a homotopy functor;
 \item $F$ stable;
 \item $F$ is half-exact.
\end{tightlist}
Then there is a unique functor $\widehat F:$ {\bf RE}
$\rightarrow$ {\bf A} such that ${\widehat F}\circ C = F$.
\end{thm}

\Pf As in the analogous theorem for $\RKK$, the only question is how
$\widehat F$ acts on morphisms.
Let $\Psi =\{\psi_t\}_{t \in [1,\infty)} : A \to B$ be an asymptotic
$C_0(X)$-morphism, and let $0 \rightarrow C_0([1,\infty), B) \rightarrow D
\stackrel{\pi}{\rightarrow} A \rightarrow 0$ be the corresponding
short exact sequence described above.  Since $C_0([1,\infty), B)$
is contractible, it follows from the long exact sequence for $F$ that
$F\pi : F(D) \rightarrow F(A)$ is an isomorphism.
Set ${\widehat F}(\Psi) = F(\ev_1)\circ\left(F(\pi)\right)^{-1}$.
It is straightforward to show that this map is well-defined,
so ${\widehat F}(\Psi)$ is a morphism
from $F(A)$ to $F(B)$.  Next, since $F$ satisfies Bott
periodicity, we have isomorphisms
among $F(A)$, $F(S^2A), F(A \otimes\K)$ and $F(S^2A\otimes\K)$ for
all $A$ in {\bf SC*(X)}, whence each element of $\RE(X;A,B)$ determines
a morphism from $F(A)$ to $F(B)$.
\EPf

\begin{thm}\label{thm:RKK-REtheory}(Compare Theorem 25.6.3, \cite{blackadar98})
 Let $A$ be a separable $C_0(X)$-algebra for which
$\RKK(X; A, B)$ is half-exact.  Then $\RE(X; A, B)$ is naturally
isomorphic to $\RKK(X; A, B)$ for every separable $C_0(X)$-algebra $B$.
In particular, if $A$ is $\RKK(X)$-nuclear in the sense of
\cite{bauval98}, then $\RE(X ; A, B) \cong \RKK(X; A, B)$ for every
separable $C_0(X)$-algebra $B$.
\end{thm}

\Pf From Theorem \ref{thm:RKKcat}, we have a homomorphism
$\alpha$ from
$\RKK(X; A, B)$ to $\RE(X ; A, B)$.  Theorem \ref{thm:REcat} implies
that we have a map $\RKK(X; A, A) \times \RE(X; A, B)
\longrightarrow \RKK(X; A, B)$, and so pairing with the element
${\bf 1}_A \in \RKK(X; A, A)$ defines a homomorphism
$\beta: \RE(X ; A, B) \longrightarrow \RKK(X; A, B)$ that is the
inverse to $\alpha$.
\EPf

\section{Examples and Applications}

In this section, we consider examples and applications of the
preceding theory. First, we consider $\RE$-elements associated to unbounded 
$\RKK$-elements, which can also arise from
families of elliptic differential operators parametrized by $X$.
We also construct \lq\lq fundamental classes\rq\rq\
for unital $C_0(X)$-algebras, and define invariants of central bimodules
in noncommutative geometry.

\subsection{Unbounded $\RKK$-elements}

Let $X$ be a locally compact space, and let $A$ and $B$ be separable
$C_0(X)$-algebras. Given any Hilbert $B$-module $\E$ with
$B$-valued inner product $\langle \ ,\ \rangle$, let $\B(\E)$
denote the $\cs$-algebra of bounded $B$-linear operators
on $\E$ that possess an adjoint. The ideal $\K(\E)$ of compact operators on
$\E$
has a natural $C_0(X)$-action \cite{kasparov88}; the structural
homomorphism $\Phi : C_0(X) \to M(\K(\E)) = \B(\E)$ is defined as
follows: For all $e \in \E$ and $f \in C_0(X)$,
$$\Phi(f)e = \lim_\lambda e \cdot (f \cdot b_\lambda),$$
where $\{b_\lambda\}$ is any approximate unit for $B$.
Thus, $\B(E)$ is a central Banach $C_0(X)$-module.

\begin{dfn}
An unbounded $\RKK(X;A,B)$-cycle is a triple $(\E, \phi, D)$ where
\begin{tightlist}
\item $\E$ is a countably generated $\Z_2$-graded Hilbert $B$-module;
\item $\psi : A \to \B(\E)$ is a $C_0(X)$-morphism;
\item $D$ is a self-adjoint regular operator on $\E$ of grading degree one
satisfying
\end{tightlist}
\begin{ilist}
\item $(D \pm i)^{-1} \psi(a) \in \K(\E)$ for all $a \in A$;
\item $\{a \in A : [D, \psi(a)] \text{ is densely defined and extends to }
\B(\E)\}$ is dense in $A$.
\end{ilist}
We define $\Psi(X; A, B)$ to be the collection of all unbounded
$\RKK(X;A,B)$-cycles.
\end{dfn}

Let $\eps$ be the grading operator of $\E$.
For any $x \in \R$ and $t \geq 1$, the operator $x \eps + t^{-1} D$ is a
self-adjoint
regular operator \cite{lance95} on $\E$ that satisfies
$$(x \eps + t^{-1} D)^2 = x^2 1_\E + t^{-2}D^2 \geq x^2 1_\E.$$
Define a one-parameter family of maps $C_0(\R) \odot A \to C_0(\R) \otimes
\K(\E)$
on elementary tensors $g \otimes a$ by the formula
$$g \otimes a \mapsto g(x \eps + t^{-1}D) \circ \psi(a)$$
and extend linearly,
where $ g(x \eps + t^{-1}D)$ is computed using the functional calculus for
self-adjoint
regular operators \cite{lance95,baaj-julg83}.

\begin{prop}
Let $(\E, \psi, D) \in \Psi(X; A, B)$. The above formula extends to
an asymptotic $C_0(X)$-morphism
$$\{\psi_t^D\} : C_0(\R) \otimes A \to C_0(\R) \otimes \K(\E),$$
which determines a well-defined map $\Psi(X; A, B) \to \RE(X; A, B).$
\end{prop}

\Pf
That this defines an asymptotic morphism is well-known
\cite{connes-higson89,higson-kasparov-trout98}; we need only
show that it is asymptotically
$C_0(X)$-linear. By an approximation argument, we only need to check
on elementary tensors $g \otimes a$. Since $\psi$ is $C_0(X)$-linear,
we have for each $f \in C_0(X)$ that
$$g(x \eps + t^{-1}D) \circ \psi(f \cdot a) = \Phi(f) g(x \eps + t^{-1}D) \circ \psi(a),$$
since $\Phi(f) \in ZM(\K(\E)) = Z(\B(\E))$. The result now easily follows.
\EPf

As in \cite{baaj-julg83} we have a map $\Psi(X; A, B) \to \RKK(X; A, B)$ given
by the formula 
$$(\E, \psi, D) \mapsto (\E, \psi, F(D)),$$
where $F(D) = D(D^2 + 1)^{-1/2}$.
Take $[(\E, \psi, T)] \in \RKK(X;A,B)$ (Definition \ref{def:RKKcycle}). We may
and do assume that $T = T^*$ and $\|T\| \leq 1$. Then
$(\E, \psi, G(T))$, where $G(x) = x(1-x^2)^{-1/2}$, is an unbounded
$\RKK(X;A,B)$-cycle which maps to
$(\E, \psi, T)$, since $(F \circ G)(x) = x$. Therefore
the map  $\Psi(X; A, B) \to \RKK(X; A, B)$ is surjective.

\begin{prop}
Let $\RKK(X; A, B) \to \RE(X; A, B)$ be the natural transformation from
Theorem \ref{thm:RKK-REtheory}. The following diagram commutes:
$$\diagram
 & \Psi(X; A, B) \dlto \drto & \\
\RKK(X; A, B) \rrto & & \RE(X; A, B).
\enddiagram$$
\end{prop}

\subsection{Elliptic Differential $C_0(X)$-operators}

For the material and definitions in this subsection, we refer the reader to
\cite{mishchenko-fomenko80} and \cite{trout99}.

Let $B$ be a separable $C_0(X)$-algebra with unit and let $M$ be a smooth
compact oriented Riemannian manifold. Let
$E \to M$ be a smooth vector $B$-bundle, i.e., a smooth locally trivial fiber
bundle with fibers $E_p \cong \E$, a finite
projective (right) $B$-module. Equip the fibers $E_p$ with smoothly-varying
$B$-valued metrics
$$\langle \cdot ,\cdot \rangle_p : E_p \times E_p \to B;$$
this is always possible since $\E$ admits partitions of unity
\cite{blackadar98}.
Let $C^\infty(E)$ denote the module of smooth sections of $E$.
We let $\H_E = L^2(E)$ denote the Hilbert $B$-module
completion of $C^\infty(E)$ with respect to the $B$-valued Hermitian metric
$$\langle s, s' \rangle = \int_M \langle s(p), s'(p) \rangle_p
\dvol_M(p),$$
where $s, s' \in C^\infty(E)$ and $\dvol_M$ is the Riemannian volume measure.

Let $D : C^\infty(E) \to C^\infty(E)$ be a self-adjoint elliptic partial
differential $B$-operator.
For simplicity, we assume that $D$ has order one and that $E$ splits as a
direct sum $E = E_0 \oplus E_1$ with respect to
a grading bundle automorphism $\epsilon$. Also, assume that $D$ is of
degree one
with respect to this grading, i.e.,
$$D = \pmatrix 0 & D^- \\ D^+ & 0 \endpmatrix$$
where $(D^-)^* = D^+$.

In the special case where $X$ is compact and $B = C(X)$,
then we can view $D = \{D_x\}_{x \in X}$ as a family of elliptic
differential operators on $M$ parametrized by $X$; this is
accomplished by letting $D_x = D \otimes_{\ev_x} 1$ on the complex
vector bundle $E^x = E \otimes_{ex} \C$, where
$\ev_x : C(X) \to \C$ denotes evaluation at $x$
(note that $C^\infty(E^x) = C^\infty(E) \otimes_{\ev_x} \C$).
Thus, for a general unital $C_0(X)$-algebra $B$, we think of $D$ as
a generalized family of elliptic operators parametrized by $X$.

Since $B$ is a $C_0(X)$-algebra and $\H_E$ is a Hilbert $B$-module, the
$\cs$-algebra of compact operators $\K(\H_E)$ is also a $C_0(X)$-algebra.
For each $x \in \R$ and $t \geq 1$, the operator $x \epsilon + t^{-1}D$ is a
self-adjoint regular operator on $\H_E$ with compact resolvents
$(D \pm i)^{-1}$.
There is a continuous family of $*$-homomorphisms \cite{trout99}
$$ \{\phi_t^D\} : C_0(\R) \to C_0(\R) \otimes \K(\H_E) :  g  \mapsto g(x
\epsilon + t^{-1}D)$$
which is defined by the functional calculus for self-adjoint regular operators
\cite{lance95}.
Using the universal property of the maximal tensor product,
we obtain an asymptotic $C_0(X)$-morphism
$$\{\hat{\phi}_t^D\} : C_0(X) \otimes C_0(\R) \to C_0(\R) \otimes \K(\H_E)$$
by mapping $f \otimes g$ to $f \cdot  g(x \epsilon + t^{-1}D)$.

\begin{dfn} We define $\lhcl D \rxhcl = \lhcl \hat{\phi}_t^D \rxhcl \in
 \RE(X; C_0(X), \K(\H_E))$.
\end{dfn}

Let $\ell^2_B$ denote the standard Hilbert $B$-module. Using the
isomorphism $\H_E \oplus \ell^2_B \cong \ell^2_B$ that comes from the
Kasparov Stabilization Theorem, we
have an inclusion $i: \K(\H_E) \hookrightarrow \K(\ell^2_B) \cong \K
\otimes B$,
where $\K$ is the algebra of compact operators on separable Hilbert space.
Thus, we can push forward the $\RE$-theory class of $D$ to obtain $i_*\lhcl D
\rxhcl \in \RE(X; C_0(X), B)$.

As an example, let $M$ be a smooth compact even-dimensional spin manifold with
Dirac operator $\DD : C^\infty(S) \to C^\infty(S)$, where
$S = S^+ \oplus S^- \to M$ is the (graded) spinor bundle. Let $F \to M$ be a
smooth vector $B$-bundle with connection
$\nabla^F : C^\infty(F) \to C^\infty(T^*M \otimes F)$. By
twisting the Dirac operator $\DD$ with the connection $\nabla^F$, we obtain
an elliptic differential $B$-operator $\DD_F : C^\infty(S \otimes F) \to
C^\infty(S \otimes F)$ associated to the vector $B$-bundle
$E = S \otimes F \to M$ as above.

\subsection{Fundamental classes of unital $C_0(X)$-algebras}

Let $A$ be a unital separable $C_0(X)$-algebra. Since $A$ is unital,  $ZM(A) =
Z(A) \subset A$, and
so we can consider the structural homomorphism
$\Phi_A : C_0(X) \to ZM(A) \hookrightarrow A$ as a canonically given
$C_0(X)$-morphism.

\begin{dfn}
Let $A$ be a unital separable $C_0(X)$-algebra. We define the {\bf fundamental
class of} $A$ to be
the $\RE$-theory element
$$\lhcl A \rxhcl = \lhcl \Phi_A \rxhcl \in \RE(X; C_0(X), A)$$
determined by the structural homomorphism.
\end{dfn}

\begin{prop}
Let $A$ and $B$ be unital $C_0(X)$-algebras. Under the $C_0(X)$-tensor product
operation
$$\RE(X; C_0(X), A) \otimes \RE(X; C_0(X), B) \to \RE(X; C_0(X), A \xotimes
B)$$
we have $\lhcl A \rxhcl \otimes \lhcl B \rxhcl = \lhcl A \xotimes B \rxhcl.$
\end{prop}

\Pf
This follows from Definition \ref{def:xtensor} and the isomorphism $C_0(X)
\xotimes C_0(X) \cong C_0(X)$.
\EPf

\begin{prop}
Let $p : Y \to X$ be a continuous map of compact spaces. Let $A$ be a unital
$C(X)$-algebra.
Under the pullback transformation
$$p^* : \RE(X; C(X), A) \to \RE(Y; C(Y), p^*A)$$
we have $p^* \lhcl A \rxhcl = \lhcl p^*A \rhcl_Y$.
\end{prop}

\Pf
Since $X$ and $Y$ are compact and $A$ is unital, $p^*A = C(Y) \xotimes
A$ is also unital.
Let $\id_Y : C(Y) \to C(Y)$ denote the identity map. From the definition of the
pullback,
$$p^*(\Phi_A) = \id_Y \xotimes \Phi_A.$$
Take $g \in C(Y)$ and $f \in C(X)$. On elementary tensors $h \xotimes a \in
p^*A$ we have
$$(\id_Y \xotimes \Phi_A)(g \xotimes f)(h \xotimes a) = gh \xotimes fa = g(f
\circ p)h \xotimes a = \Phi_{p^*A}(g(f \circ p))(h \xotimes a).$$
Since $p^*C(X) \cong C(Y)$,  we see that $p^*(\Phi_A) = \Phi_{p^*A}$,
as desired.
\EPf

\subsection{Central Bimodules}

Let $A$ be a separable $\cs$-algebra with unit and center $Z(A)$ (we
restrict to unital $\cs$-algebras so that we are guaranteed that
$A$ has a nonempty center).

\begin{dfn}
A {\bf central bimodule} over $A$ is an $A$-$A$-bimodule $\E$ such that
$z \cdot e = e \cdot z$ for all $e \in \E$ and $z \in Z(A)$.
\end{dfn}

Central bimodules have found several applications in noncommutative
geometry and
noncommutative physics \cite{d-v-michor96,d-h-l-s96,masson96}.
This is due to the fact that the module $\E = \Gamma(E)$ of continuous sections
of a complex vector bundle $E \to X$
on a compact space $X$ naturally defines a bimodule over the $\cs$-algebra $A =
C(X)$ for which the left and right
actions are compatible. Allowing $A$ to be noncommutative, we see that the
notion of a central bimodule is a
slight generalization of the notion of a complex vector bundle to the
noncommutative setting.

Let $X = \widehat{Z(A)}$ denote the spectrum of the center of $A$, i.e., $Z(A)
\cong C(X)$ (note that $X$
is compact since $Z(A)$ is unital). This gives
$A$ the structure of a $C(X)$-algebra. Thus, if $(\E, \psi, T)$ is a
(unbounded) $\RKK$-cycle, then
$\E$ is naturally endowed with the structure of a central bimodule over $A$.

Dubois-Violette and Michor \cite{d-v-michor96} have asked how to
define a suitable notion of $K$-theory for $A$ which
is appropriate for defining invariants for central bimodules. The previous
observations lead us to propose the following definition.

\begin{dfn}
Let $A$ be a $\cs$-algebra with unit. We define the {\bf central $K$-theory}
$ZK(A)$ of $A$ to be
$$ZK(A) =  \RE(\widehat{Z(A)}; A, A).$$
If $A$ is $\RKK(X)$-nuclear (see Theorem \ref{thm:RKK-REtheory}),
then $ZK(A) \cong \RKK(\widehat{Z(A)}; A, A)$.
\end{dfn}

The relationship between $K$-theory and central $K$-theory is contained in the
following proposition, which gives the two extremes for the center of $A$.

\begin{lem}
Let $A$ be a $\cs$-algebra with unit. Then $ZK(A)$ is a ring with
multiplication
given by
the composition product. If $A = C(X)$ is commutative, then $ZK(A) \cong
K_0(A)$. If $Z(A) = \C$ and $A$ is
$K$-nuclear, then $ZK(A) \cong KK(A,A)$.
\end{lem}

\Pf That $ZK(A)$ is a ring with multiplication given by the composition product
follows from Theorem \ref{thm:REtheory}. If $A \cong C(X)$ then
$ZK(A) = \RE(X; C(X), C(X)) \cong \RKK(X; C(X), C(X)) \cong K^0(X) \cong
K_0(A)$
by Proposition 2.20 in \cite{kasparov88}. If $Z(A) = \C$ then
$X = \widehat{Z(A)} = \bullet$ and so $ZK(A) = \RE(\bullet, A, A) = E(A,A)$ and
the result follows. \EPf

Note that the ring structure incorporates the fact that if $\E$ and $\F$ are
central bimodules over $A$ then
$\E \otimes_A \F$ is also a central bimodule over $A$.


\bibliographystyle{amsalpha}
\providecommand{\bysame}{\leavevmode\hbox to3em{\hrulefill}\thinspace}

\end{document}